\documentclass[11pt,twoside]{article}
\usepackage{amsfonts,amsmath,epsfig,amssymb,amsthm}

\usepackage{cite}

\usepackage{xcolor}

\numberwithin{equation}{section}
\newfont{\ctv}{msam10}
\newcommand{\bbox}{\mbox{\ctv \symbol{3}}}
\def\QED{{$\hfill\bbox$}}
\newenvironment{pf}[1]{\par\vskip1mm{\noindent\it #1.}\ }{\QED\par\vskip2mm}
\newtheorem{theorem}{Theorem}[section]

\newtheorem{lemma}{Lemma}[section]
\newtheorem{definition}{Definition}[section]

\newtheorem*{theoremH}{Hypotheses (H)}
\newtheorem*{theoremU}{Hypotheses (U)}
\newtheorem*{theoremq}{Hypotheses (q)}

\newcommand{\Barint}{-\kern-.15in\int}
\newcommand{\barint}{-\kern-.12in\int}
\renewcommand{\r}{\mathbb{R}}
 
\DeclareMathOperator{\haus}{haus} 
\DeclareMathOperator{\co}{co} 
\DeclareMathOperator{\dom}{dom}

\def\bpf{\begin{pf}}
\def\epf{\end{pf}}

\begin{document}

\thispagestyle{empty}
 \setcounter{page}{1}

 \title{Optimal control of a population dynamics model with hysteresis \footnote{Supported by National Natural Science Foundation of China (12071165 and 62076104), Natural Science Foundation of Fujian Province (2020J01072), Program for Innovative Research Team in Science and Technology in Fujian Province University, Quanzhou High-Level Talents Support Plan (2017ZT012), and by Scientific Research Funds of Huaqiao University (605-50Y19017, 605-50Y14040). The research of the first author was also supported by Ministry of Science and Higher Education of Russian Federation (Grant No. 075-15-2020-787, large
scientific project ``Fundamentals, methods and technologies for
digital monitoring and forecasting of the environmental situation on
the Baikal natural territory'').}
}
\author{Sergey A. Timoshin\footnote{Fujian Province University Key Laboratory of Computational Science, School of Mathematical Sciences, Huaqiao University, Quanzhou 362021, China -- and --
Matrosov Institute for System Dynamics and Control
 Theory,
 Russian Academy of Sciences, Lermontov str. 134,
664033  Irkutsk,
 Russia, E-mail: {\tt
sergey.timoshin@gmail.com}. } $^{,}$\footnote{Corresponding
author.}
\ and
Chen Bin\footnote{Fujian Province University Key Laboratory of Computational Science, School of Mathematical Sciences, Huaqiao University, Quanzhou 362021, China, E-mail: {\tt
chenbinmath@163.com}. }
}

\date{}

\maketitle \hyphenation{Ca-ra-the-odo-ry}

\begin{abstract}
\noindent {This paper addresses a nonlinear partial differential control system arising in
population dynamics. The system consist of three diffusion equations describing the evolutions of three biological species: prey, predator, and food for the prey or vegetation. The equation for the food density incorporates a hysteresis operator of generalized stop type accounting for underlying hysteresis effects occurring in the dynamical process. We study the problem of minimization of a given integral cost functional over solutions of the above system. The set-valued mapping defining the control constraint is state-dependent and its values are nonconvex as is the cost integrand as a function of the control variable. Some relaxation-type results for the minimization problem  are obtained and the existence of a nearly optimal solution is established.}
\end{abstract}

\small\textbf{Keywords}:  optimal control problem;  hysteresis; biological diffusion models; nonconvex integrands; nonconvex control constraints.


\vspace{0.3cm}

\small\textbf{2010 Mathematics Subject Classification}: 49J20, 49J21, 49J45,
49J53.

\vspace{0.3cm}

\small\textbf{Abbreviated title}: Optimal control of a biomodel with hysteresis.

\section{Introduction} \label{number of section:introduction}

Given a  final time $T>0$ and  a bounded domain $\Omega\subset
\mathbb{R}^N$, $N\leq 3$, with smooth
boundary $\partial\Omega$, the primary concern of this paper is to consider the following optimal control problem:
\begin{equation*}J(\sigma,v,w,u):=\int_0^T \int_\Omega q(t,x,\sigma(t,x),v(t,x),w(t,x),u(t,x)) \, dx dt\rightarrow \inf  \eqno (P)\end{equation*}

over the solution set of the dynamical control system:
\begin{equation}
\sigma_t -(\lambda(v))_t-\kappa \Delta\sigma + \partial I_{v,w}(\sigma) \ni F(\sigma,v,w)\,  u \qquad \text{in} \; Q(T)
,\label{1.1}
\end{equation}
\begin{equation} v_t - \Delta v=h(\sigma, v,w) \qquad \text{in}
\; Q(T) , \label{1.2}\end{equation}
\begin{equation} w_t - \Delta w=g(\sigma, v,w)  \qquad \text{in}
\; Q(T) , \label{1.3}\end{equation}
\begin{equation}
\sigma(x,0)=\sigma_0(x), \quad v(x,0)=v_0(x), \quad w(x,0)=w_0(x)   \quad \text{in} \; \Omega ,\label{1.4}
\end{equation}
\begin{equation}
\frac{\partial \sigma}{\partial n}=\frac{\partial v}{\partial n}=\frac{\partial w}{\partial n}=0 \quad \text{on} \; \partial\Omega\times [0,T] ,\label{1.5}
\end{equation}

subject to the state-dependent control constraint:
\begin{equation}
u\in U(t,x,\sigma, v,w) \qquad \text{in} \; Q(T) .\label{1.6}
\end{equation}

Here,  $Q(T) := [0,T]\times \Omega$, $\kappa$ is a given constant, $\sigma_0$, $v_0$,  $w_0$ are given  initial
conditions,  $\partial/\partial n$ is the outward normal derivative on $\partial\Omega$. The operator $\partial I_{v,w}(\cdot)$ is
the subdifferential, in the sense of convex analysis, of the
indicator function $I_{v,w}(\cdot)$  of the interval $[f_*(v,w),
f^*(v,w)]$, $f_*,f^*:\mathbb{R}^2\to \mathbb{R}$, $q:Q(T)\times\mathbb{R}^4\to \mathbb{R}$, $F, h, g:\mathbb{R}^3\to \mathbb{R}$, $\lambda:\mathbb{R}\to \mathbb{R}$ are prescribed functions,  $U:Q(T)\times \mathbb{R}^3\to \mathbb{R}$ is a multivalued mapping with compact, but not necessarily convex, values.

\hspace{0.6cm} We note that when $w$ is fixed and $F\equiv 0$, $\lambda(v)=v$ in $(\ref{1.1})$, inclusion $(\ref{1.1})$ yields the differential representation of hysteresis operator  of the generalized stop type (cf. \cite{Visintin}). During the last three decades, hysteresis operators have been extensively  applied in modeling, analysis, and control of a variety of irreversible nonlinear phenomena in applied sciences including phase transitions \cite{CKK, KSS, Giorgi, HH, CF, TCON2018}, porous media flow \cite{Kubo, KOKPR, Albers, DKR, SICON2016},  thermostat models \cite{CSSZ, KK, LRS, GR, ESAIM}, concrete carbonation \cite{AK, AK2, JJG, JMMA2017, NONRWA2019} and many others. In the same vein, system $(\ref{1.1})$--$(\ref{1.5})$ (without the control $u$) was introduced in \cite{AM} (see, also, \cite{ZW, WZ}) to model the evolution of populations in the vegetation-prey-predator framework when diffusive effects in the dynamics of three species are taken into account and the food density for the prey exhibits   hysteretic character. This latter means that the growth rate of the food for the prey depends not only on the current state of prey and predator populations, but it also depends on their immediately preceding density history. In \cite{SCL2019}, with the aim to achieve  a possible optimization of the population dynamics process by way of controlling the food supply for the prey this uncontrolled system was subjected to a control action and the existence of solutions for the corresponding control system $(\ref{1.1})$--$(\ref{1.6})$ was established.

\hspace{0.6cm} In our model, the unknown variables  $\sigma$,  $v$, and $w$ represent the densities of vegetation,  prey and  predator, respectively, with Eqs. $(\ref{1.1})$, $(\ref{1.2})$, and $(\ref{1.3})$ characterizing the evolutions of the corresponding species. In this respect, the nonconvexity of values of $U$ in the control constraint $(\ref{1.6})$ is a biologically relevant assumption.

\hspace{0.6cm} Our aim in this paper is to consider along with $(P)$ the following alternative problem:
\begin{equation*}J_U^{**}(\sigma,v,w,u):=\int_0^T \int_\Omega q_U^{**}(t,x,\sigma(t,x),v(t,x),w(t,x),u(t,x)) \, dx dt\rightarrow \inf  \eqno (RP)\end{equation*}

over  solutions of  $(\ref{1.1})$--$(\ref{1.5})$ supplemented with the following alternative control constraint:
\begin{equation} u\in \co  U(t,x,\sigma, v,w) \qquad \text{in} \; Q(T),\label{1.7}\end{equation}

and explore some properties and relationships between Problems $(P)$ and $(RP)$. Here, $(RP)$ stands for the \emph{relaxed problem}, $\co  U$ denotes  the convex hull of the set $U$, which is the intersection of all convex sets containing $U$, and for the extended real-valued
function $q_{U}:Q(T)\times
\mathbb{R}^3\times\mathbb{R}\rightarrow\mathbb{R}\cup\{+\infty\}$:
\begin{equation*}q_{U}(t,x,\sigma,v,w,u)= \left\{
\begin{array}{ll}
q(t,x,\sigma,v,w,u) & \mbox{if} \;u \in U(t,x,\sigma,v,w) ,\\
+\infty  & \mbox{otherwise},
\end{array}\right.
\label{}\end{equation*}  $q_U^{**}(t,x,\sigma,v,w,u)$ is the bipolar or the second
conjugate of the function $u\mapsto q_{U}(t,x,\sigma,v,w,u)$, which is the largest lower semicontinuous convex function less than or equal to
$q_U^{**}$. We note that both the values of the set-mapping in the control constraint and the integrand as a function of control in Problem $(RP)$ are  convex.

\hspace{0.6cm} We remark that while our model has not been inspired by any particular applied problem, the results we obtain might find potential applications in analysis and control of real-world ecosystems such as spruce budworm population dynamics models. The latter describe interactions in the budworm-forest framework and include three species: boreal forest trees such as balsam fir, insect pests such as spruce budworm which feed on the foliage of trees and avian predators feeding in turn on the insects. These ecosystems are typical to Northern parts of Canada and Russia, including the lake Baikal region, where the defoliation of larches by insects and the subsequent dieback of the trees is a major problem in forest ecology. Available control actions towards preservation of forest reserve include spraying insecticides, removal of infected trees and other and affect the rate of change in the insects population. This prompted us to include the control function $u$ into the first equation of our evolution system. On the other hand, the natural concern to avoid the total perishing of birds population caused by insufficient food supply should be reflected in optimization strategy towards minimization of forest loss.

\hspace{0.6cm} The outline of the paper is as follows. After introducing in the next section the notation and hypotheses on the data describing Problem $(P)$, in Section 3 we derive several continuity properties of the control-to-state solution operator associated with our control systems. These properties are instrumental in proving our main results in the last Section 4. The latter consist in establishing the existence of an optimal solution for Problem $(RP)$, the density, in an appropriate topology, of solutions of the original control system $(\ref{1.1})$--$(\ref{1.6})$ among solutions of the convexified control system $(\ref{1.1})$--$(\ref{1.5})$, $(\ref{1.7})$, and the existence of a nearly optimal, in a suitable sense, solution for problem $(P)$.

\hspace{0.6cm} At the end of Introduction, we mention that the existing literature on optimal control of systems exhibiting hysteretic character is sparse. It is worthwhile however to mention a seminal contribution of M. Brokate to this area (see \cite{Brokate1, Brokate2}). Moreover, to the best of the authors' knowledge, there have been no contributions so far addressing optimal control problems for biological models with hysteresis.

\vspace{0.3cm}

\textbf{Acknowledgements.}  The authors want to thank the
anonymous referees for their valuable suggestions and remarks
which helped to improve the manuscript.

\section{Preliminary notions, hypotheses on the data, and statement of the main results} \label{section2}

In this section, we fix the notation which we use throughout the paper, prove an auxiliary lemma, specify hypotheses which we posit on the data describing Problem $(P)$, give a precise meaning in which solutions to our control systems are understood and the corresponding minimization problems are treated, and state our main results.

\hspace{0.6cm}  Let $X$ be a Hilbert space with the inner product $(\cdot ,\cdot )_X$. A function $\varphi:X\rightarrow\mathbb{R\cup\{+\infty\}}$ is called proper
if its effective domain $\dom \varphi:=\{x\in X; \varphi(x)< +\infty \}$
is nonempty. By definition, the subdifferential $\partial
\varphi(x)$, $x\in X$, of a proper, convex, lower semicontinuous
function $\varphi$ is the set
$$\partial
\varphi(x)=\{h\in X; ( h, y-x)_X\leq
\varphi(y)-\varphi(x) , \; \forall y\in X\} ,$$

and its domain is the set $\dom\partial\varphi:=\{x\in X; \partial\varphi(x)\neq\emptyset \}$. It is known \cite{Brezis} that $\dom\partial\varphi\subset\dom\varphi$.

\hspace{0.6cm} We say that a sequence of proper, convex, lower
semicontinuous functions $\varphi_n:X\to \mathbb{R}, \;
n\geq 1$, \emph{Mosco-converges} \cite{Mosco} to a proper, convex,
lower semicontinuous function $\varphi:X\to \mathbb{R}$,
denoted $\varphi_n\stackrel{M}{\to}\varphi$, if:
\begin{itemize}
\item[1)] for any $x\in X$ and any sequence $x_n\in X, \; n\geq 1$,
weakly converging to $x$ we have
$$
\varphi (x)\leq \liminf\limits_{n\to\infty}\varphi_n(x_n);
$$
\item[2)] for any $x\in X$ there exists a  sequence $x_n\to x$
such that $\varphi_n(x_n)\to \varphi (x)$.
\end{itemize}

\hspace{0.6cm} Denote by $d_X(x,A)$ the distance from a point $x\in X$ to a set
$A\subset X$. Then, the Hausdorff metric  on the space of closed
bounded subsets of $X$, denoted $cb (X)$, is the function:
$$\haus_X(A,B)=\max\{\sup\limits_{x\in A}d_X(x,B),\sup\limits_{y\in B}d_X(y,A)\},
\quad A, B \in cb(X).$$

\hspace{0.6cm} If $({\mathcal E},{\mathcal A})$ is a measurable
space, then a multivalued mapping $F: {\mathcal E}\to
cb(X)$ is called measurable if $\{\tau \in {\mathcal E}; \;
F(\tau )\cap C\neq \emptyset \}\in {\mathcal A}$ for any closed
subset $C$ of $X$. A set $\mathcal{F}$ of measurable
functions from $\mathcal{E}$ to $X$ is called decomposable if for
any $f_1,f_2\in \mathcal{F}$ and any $E\in \mathcal{A}$ we have that
$f_1\cdot \chi_E +f_2\cdot \chi_{\mathcal{E}\setminus E}\in
\mathcal{F}$, where $\chi_E$ stands for the characteristic function
of the set $E$.

\hspace{0.6cm} Let $Y, Z$ be two Banach spaces. A multivalued mapping $G:Y\to Z$ is called lower semicontinuous if for any $y\in Y$, $z\in F(y)$ and any sequence $y_k\to y$, $k\geq 1$, there exists a sequence
$z_k\in F(y_k)$, $k\geq 1$, which converges to $z$.

\hspace{0.6cm} Given two functions $f_*, f^*:\mathbb{R}^2\to \mathbb{R}$, recall that the
indicator function $I_{v,w}(\cdot)$ of the
set $K(v,w):=[f_{*}(v,w), f^{*}(v,w)]$ is defined as follows
\begin{equation*} I_{v,w}(\sigma):=  \left\{
\begin{array}{cl}
0 &\mbox{if} \;\; \;  \sigma\in K(v,w) ,  \\
 +\infty &\mbox{otherwise} .
\end{array}
\right.
\end{equation*}

Its subdifferential has the form:

\begin{equation*}\partial I_{v,w}(\sigma)= \left\{ \begin{array}{ccl}
\emptyset &\mbox{if}& \sigma \notin K(v,w) ,\\
{[}0,+\infty) &\mbox{if}& \sigma = f^{\ast}(v,w)> f_{\ast}(v,w),\\
\{0\}       &\mbox{if}& f_{\ast}(v,w)< \sigma < f^{\ast}(v,w) ,\\
(-\infty, 0{]} &\mbox{if}& \sigma = f_{\ast}(v,w)< f^{\ast}(v,w),\\
(-\infty, +\infty) &\mbox{if}& \sigma = f_{\ast}(v,w)= f^{\ast}(v,w).
\end{array}
\right. \label{}
\end{equation*}

And, for $\mu>0$, the Yosida regularization of $\partial I_{v,w}(\sigma)$ is
the function
\begin{equation*}\partial I^\mu_{v,w}(\sigma)= \frac{1}{\mu}[\sigma-f^*(v,w)]^+-\frac{1}{\mu}[f_*(v,w)-\sigma]^+ , \quad \sigma,v,w\in \mathbb{R}.\label{}
\end{equation*}

\hspace{0.6cm} Let $H$ be a  Hilbert space  with the  norm $|\cdot|_H$.

\begin{lemma} \label{lemma2.1} Let $v_n\to v$, $w_n\to w$, $\sigma_n\to \sigma$ in $L^2(0,T;H)$, $f_n\to f$ weakly in
$L^2(0,T;H)$. In addition, assume that a set-valued function $\mathcal{K}:H\times H\to cb (H)$ has convex values and
\begin{equation}\haus_H({\mathcal{K}}(v_1,w_1),{\mathcal{K}}(v_2,w_2))  \leq
    R\, (|v_1-v_2|_H+|w_1-w_2|_H),
    \label{2.1}\end{equation}

$v_i,w_i\in H$, $i=1,2$, where $R>0$ is a constant. If
\begin{equation*}f_n(t)\in\partial I_{\mathcal{K}(v_n(t),w_n(t))}(\sigma_n(t)) \hspace{0.3cm} \text{for a.e. }\; t\in [0,T] , \label{}\end{equation*}  then
\begin{equation*}f(t)\in\partial
I_{\mathcal{K}(v(t),w(t))}(\sigma(t))\hspace{0.3cm} \text{for a.e. }\; t\in [0,T]
,\label{}\end{equation*}
where $I_{\mathcal{K}(v(t),w(t))}$ is the indicator function of the set ${\mathcal{K}(v(t),w(t))}\subset H$, $t\in [0,T]$.
\end{lemma}

\bpf{Proof} First, we define on $L^2(0,T;H)$ the function
$$I_{\mathcal{K}(v,w)}(\sigma)=\int_0^T I_{\mathcal{K}(v(\tau),w(\tau))}(\sigma(\tau))\, d\tau.$$

It is easy to see that $I_{\mathcal{K}(v,w)}$ is the indicator function of the set
$${\mathcal{K}(v,w)}=\{\sigma\in L^2(0,T;H); \; \sigma(t)\in \mathcal{K}(v(t),w(t)) \hspace{0.2cm} \text{a.e.
on}\;[0,T]\}$$

and it is proper, convex, and lower semicontinuous. Moreover, \cite[Proposition 0.3.3]{Kenmochi} implies that
\begin{equation}f_n\in\partial I_{\mathcal{K}(v_n,w_n)}(\sigma_n), \; \; n\geq 1 . \label{2.2}\end{equation}

\hspace{0.6cm}  Next, we prove that  $I_{\mathcal{K}(v_n,w_n)}\stackrel{M}{\to}I_{\mathcal{K}(v,w)}$. According to the definition of the Mosco-convergence introduced above, we need to show that
for  any sequence $z_n\to z$
weakly in  $L^2(0,T;H)$  we have
\begin{equation}
I_{\mathcal{K}(v,w)}(z)\leq \liminf\limits_{n\to\infty}I_{\mathcal{K}(v_n,w_n)}(z_n),
\label{2.3}\end{equation}

and for any $z\in L^2(0,T;H)$ there exists a  sequence $z_n\to z$ in $L^2(0,T;H)$
such that
\begin{equation}
I_{\mathcal{K}(v_n,w_n)}(z_n) \to I_{\mathcal{K}(v,w)}(z).
\label{2.4}\end{equation}

\hspace{0.6cm}  So, let $z_n\to z$
weakly in  $L^2(0,T;H)$. In the case when $\liminf\limits_{n\to\infty}I_{\mathcal{K}(v_n,w_n)}(z_n)=+\infty$, (\ref{2.3}) trivially holds. Hence, without loss of generality, we can assume that  $I_{\mathcal{K}(v_n,w_n)}(z_n)=0$, $n\geq 1$, which implies that $z_n\in \mathcal{K}(v_n,w_n)$. Invoking the Mazur lemma from this inclusion and (\ref{2.1}) we infer that
$$ z\in \bigcap\limits_{k=1}^\infty \overline{\co} \left(\bigcup\limits_{n=k}^\infty z_n\right)\subset \overline{\co} \, \mathcal{K}(v,w) =\mathcal{K}(v,w).$$

Therefore, $I_{\mathcal{K}(v,w)}(z)=0$ and (\ref{2.3}) follows.

\hspace{0.6cm} Take now an arbitrary $z\in L^2(0,T;H)$. When $I_{\mathcal{K}(v,w)}(z)=+\infty$, (\ref{2.3}) implies that (\ref{2.4}) holds for any sequence $z_n\to z$ in $L^2(0,T;H)$. Consequently, we assume that
$I_{\mathcal{K}(v,w)}(z)=0$, i.e. $z\in \mathcal{K}(v,w)$. From (\ref{2.1}) it follows that we can find a sequence $z_n\to z$ in $L^2(0,T;H)$,  $z_n\in \mathcal{K}(v_n,w_n)$. In this case, $I_{\mathcal{K}(v_n,w_n)}(z_n)=0$ and (\ref{2.4}) holds again.

\hspace{0.6cm} From (\ref{2.2}) we deduce that $\sigma_n\in \dom I_{\mathcal{K}(v_n,w_n)}={\mathcal{K}(v_n,w_n)}$, i.e.
$I_{\mathcal{K}(v_n,w_n)}(\sigma_n)=0$. From the Mosco-convergence established above we obtain
$$0\leq I_{\mathcal{K}(v,w)}(\sigma)\leq \liminf\limits_{n\to\infty}I_{\mathcal{K}(v_n,w_n)}(\sigma_n)=0 .$$

Hence, $\sigma\in \dom I_{\mathcal{K}(v,w)}$. Furthermore, for any $z\in \dom
I_{\mathcal{K}(v,w)}$ there exists a sequence $z_n\in \dom I_{\mathcal{K}(v_n,w_n)}$, $z_n\to
z$ in $L^2(0,T;H)$ such that $I_{\mathcal{K}(v_n,w_n)}(z_n)\to
I_{\mathcal{K}(v,w)}(z)$. The definition of the subdifferential together with (\ref{2.2})  imply that
$$( f_n,z_n-\sigma_n)_{L^2(0,T;H)}\leq I_{\mathcal{K}(v_n,w_n)}(z_n)-I_{\mathcal{K}(v_n,w_n)}(\sigma_n)=0 .$$
Passing to the limit in this inequality  we obtain
\begin{equation*}( f,z-\sigma)_{L^2(0,T;H)}\leq 0=I_{\mathcal{K}(v,w)}(z)-I_{\mathcal{K}(v,w)}(\sigma). \label{}\end{equation*}

Since  $z\in \dom I_{\mathcal{K}(v,w)}$ is arbitrary, from
this inequality  it follows that
$$f\in\partial I_{\mathcal{K}(v,w)}(\sigma) .$$

The claim of the lemma finally follows from \cite[Proposition 0.3.3]{Kenmochi}. \epf

\hspace{0.6cm} In the rest of the paper,  $H$ denotes the Hilbert
space $L^2(\Omega)$ endowed with the standard inner product $(\cdot,\cdot)_H$
and the associated norm $|\cdot|_H$, and  $V$ denotes the Sobolev space $H^1(\Omega)$ equipped with the norm
$|v|_V=( v,v)_V^{1/2}$,
where $( v,w)_V=( v,w)_H+a(v,w),$
\begin{equation*}
a(v,w)=\int_{\Omega}  \nabla v(x)\cdot \nabla w(x)\,
dx, \; \; v,w\in V,
\end{equation*}

$H^2(\Omega)$ is the Sobolev space $W^{2,2}(\Omega).$  Consider the linear continuous operator ${\mathcal
L}:V\to V'$ defined by
\begin{equation*}
 \langle {\mathcal L}v,w\rangle =a(v,w), \; \; v,w\in V,
\end{equation*}
where $V'$ is the dual space of $V$ and $\langle \cdot ,\cdot \rangle$ is the bilinear form
establishing the duality between $V$ and $V'$. Let $-\Delta_N: \; D(-\Delta_N)\subset H\to H$
denote the restriction of the operator ${\mathcal L}$ to the
set of elements $v\in V$ such that ${\mathcal L}v \in H$. Then,
$
D(-\Delta_N)= \left\{v\in H^2(\Omega); \; {\partial v}/{\partial n}
= 0 \; \; \hbox{in} \;\; H^{1/2}(\partial \Omega ) \right\}$
and
$-\Delta_N v=-\Delta v$ for all $v\in D(-\Delta_N)$.

\hspace{0.6cm} Problem $(P)$ is considered under the following hypotheses:

\begin{theoremH} The following assumptions hold throughout the paper:
\begin{itemize}
\item[\textup{\textbf{(H1)}}]  $\kappa>0$ is a given
    constant,  $\lambda\in C^2(\mathbb{R})$ is a given
    function with $\lambda'$, $\lambda''$ bounded on
    $\mathbb{R}$;

\item[\textup{\textbf{(H2)}}] $f_*, f^*\in C^2(\mathbb{R}^2)\cap
W^{2,\infty}(\mathbb{R}^2)$ with $0\leq f_*\leq f^*\leq 1$ on
$\mathbb{R}^2$, and $h(\sigma,0,w)=0$ for $\sigma\in [0,1]$, $w\in \mathbb{R}$, $g(\sigma,v,0)=0$ for $\sigma\in [0,1]$, $v\in \mathbb{R}$;

\item[\textup{\textbf{(H3)}}] $F, h, g:\mathbb{R}^3\to \mathbb{R}$ are locally
Lipschitz continuous functions;

\item[\textup{\textbf{(H4)}}] $\sigma_0, v_0, w_0\in L^{\infty}(\Omega)\cap V$ with $v_0\geq
0$, $w_0\geq 0$ and $f_*(v_0,w_0)\leq \sigma_0\leq f^*(v_0,w_0)$ a.e. on $\Omega$.

\end{itemize}
\end{theoremH}

\hspace{0.6cm} We note that the bounds for $f_*$ and $f^*$ in $(H2)$ are justified from the biological point of view. Indeed, when the prey population $v$ is zero,
the vegetation $\sigma$ in our three-species model stays constant, say it is one, after rescaling. And if $v$ reaches some excessive number, then all the vegetation is devoured. So, we may assume that
$\sigma=0$ in this case.

\hspace{0.6cm} Denote by ${{\mathbb R}}^+:=[0,+\infty)$. In connection with the constraint $(\ref{1.6})$, we assume the
following:

\begin{theoremU} The multivalued mapping $U:[0,T]\times \Omega \times {{\mathbb R}}^3
\to {cb}({\mathbb R})$ has the following properties:

\begin{itemize}
\item[\textup{\textbf{(U1)}}] the mapping $(t,x)\to U(t,x,\sigma,v,w), \; \sigma,v,w\in {\mathbb R},$
is measurable;

\item[\textup{\textbf{(U2)}}]  there exists  $k\in
L^2(0,T;{{\mathbb R}}^+)$ such that
\begin{align*} \label{}
 \haus_{\mathbb R}(U(t,x,\sigma_1,v_1,w_1)&,U(t,x,\sigma_2,v_2,w_2)) \\ &\leq k(t) (|\sigma_1-\sigma_2|+|v_1-v_2|+|w_1-w_2|)
\end{align*}

a.e. on $Q(T)$, $\sigma_i,v_i,w_i\in \r$, $i=1,2$;

\item[\textup{\textbf{(U3)}}] there exists a constant $m>0$ such that
\begin{equation*} \label{}
  |U(t,x,\sigma,v,w)| \leq m \quad
\mbox{a.e. on} \;\, Q(T), \sigma,v,w \in \mathbb{R}.
\end{equation*}
\end{itemize}
\end{theoremU}

The last set of hypotheses lists the assumptions we impose on the cost integrand:

\begin{theoremq} The function $q:\Omega_T\times \mathbb{R}^3\times \mathbb{R}\rightarrow\mathbb{R}$ is such that

\begin{itemize}
\item[\textup{\textbf{(q1)}}] the function $(t,x)\mapsto q(t,x,\sigma,v,w,u)$ is
measurable;

\item[\textup{\textbf{(q2)}}] there exist $k_0\in
L^2(Q(T);{{\mathbb R}}^+)$  such that
\begin{align*} \label{}
 |q(t,x,\sigma_1,v_1,w_1,u_1)&-q(t,x,\sigma_2,v_2,w_2,u_2)| \\ &\leq k_0(t,x)
 (|\sigma_1-\sigma_2|+|v_1-v_2|+|w_1-w_2|+|u_1-u_2|)
\end{align*}
a.e. on $Q(T)$, $\sigma_i,v_i,w_i\in \r$, $|u_i|\leq m$, $i=1,2;$

\item[\textup{\textbf{(q3)}}] there exist functions $k_i\in L^1(Q(T);\mathbb{R}^+)$, $i=1,2$, $k_3\in L^2(Q(T);\mathbb{R}^+)$ such that
 $$|q(t,x,\sigma,v,w,u)|\leq k_1(t,x)+k_2(t,x) |u|+k_3(t,x) (|\sigma|+|v|+|w|)$$
a.e. on $\Omega_T$, $\sigma,v,w\in \mathbb{R}$, $|u|\leq m$.

\end{itemize}
\end{theoremq}

\hspace{0.6cm} Next, we reformulate our problem
$(P)$ in a function spaces framework. To this end,  define the multivalued mapping
\begin{equation*}
\label{} {\mathcal U}(t,\sigma,v,w)=\{u\in H; \; u(x)\in U(t,x,\sigma(x),v(x),w(x))
\; \; \hbox{a.e. on} \; \; \Omega \}, \quad  \sigma, v, w \in H,
\end{equation*}
the function
\begin{equation*} \hspace{-2.1cm}
\label{} \mathfrak{q}(t,\sigma,v,w,u)=\int_\Omega
q(t,x,\sigma(x),v(x),w(x),u(x))\, dx, \quad  \sigma, v, w, u \in H,
\end{equation*}
and the set
$${\mathcal{K}}(v,w)=\{\sigma\in H; \; f_*(v(x),w(x))\leq \sigma(x) \leq f^*(v(x),w(x)) \;\;\;\text{a.e. on}\; \Omega\} , \quad v,w\in H,$$

Theorem 1.5 in \cite{HiaiUmegaki} implies that
\begin{equation*}
\label{} \overline{\rm co}\,{\mathcal U}(t,\sigma,v,w)=\{u\in H; \;
u(x)\in {\rm co}\,U(t,x,\sigma(x),v(x), w(x)) \; \; \hbox{a.e. on} \; \;
\Omega \}.
\end{equation*}

Given Hypotheses $(U)$ and $(q)$ it is a routine matter to
verify (cf., e.g., \cite[Lemmas 3.1 and 3.2]{TCON2018}) that the  mapping
${\mathcal{U}}:[0,T]\times H^3\to {cb}(H)$ has the
properties:
\begin{itemize}
\item[$\textbf{(}\mathcal{U}\textbf{1)}$] the mapping $t\mapsto
    {\mathcal{U}}(t,\sigma,v,w)$ is measurable, $\sigma,v,w\in
    H$;

\item[$\textbf{(}\mathcal{U}\textbf{2)}$] \*
\vspace{-0.8cm}\begin{align*}
\hspace{-1.3cm}    \haus_H({\mathcal{U}}(t,\sigma_1,v_1,w_1)&,{\mathcal{U}}(t,\sigma_2,v_2,w_2))  \\ &\leq
    k(t)(|\sigma_1-\sigma_2|_H+|v_1-v_2|_H+|w_1-w_2|_H) \end{align*}
    a.e. on $[0,T]$,  $\sigma_i,v_i,w_i\in H$, $i=1,2$, for $k\in L^2(0,T;\mathbb{R}^+)$ as above;

\item[$\textbf{(}\mathcal{U}\textbf{3)}$] $|{\mathcal{U}}(t,\sigma,v,w)|_{H}\leq
    m \mu(\Omega)^{1/2}$\;\;  a.e. on $[0,T]$,  $\sigma,v,w\in
    H$, where $m>0$ is as above and $\mu(\Omega)$ is the Lebesgue measure of $\Omega$,
\end{itemize}

and the function $\mathfrak{q}:[0,T]\times H^3\times
H\to \mathbb{R}$  has the properties:
\begin{itemize}
\item[$\textbf{(}\mathfrak{q}\textbf{1)}$] the function $t\mapsto
    \mathfrak{q}(t,\sigma,v,w,u)$ is measurable, $\sigma,v,w,u\in H$;

\item[$\textbf{(}\mathfrak{q}\textbf{2)}$] \*\vspace{-0.78cm}\begin{align*} \label{}
\hspace{-0.2cm}    |\mathfrak{q}(t,\sigma_1,v_1,w_1,&u_1)-\mathfrak{q}(t,\sigma_2,v_2,w_2,u_2))\\
    &\leq
    k_0^*(t)(|v_1-v_2|_{H}+|w_1-w_2|_{H}+|\sigma_1-\sigma_2|_{H}+|u_1-u_2|_{H})\end{align*}
    a.e. on $[0,T]$,  $\sigma_i,v_i,w_i,u_i\in H$,
    $i=1,2$;

\item[$\textbf{(}\mathfrak{q}\textbf{3)}$] \*\vspace{-0.43cm}\begin{equation*}
\hspace{-2.34cm}    \label{}|\mathfrak{q}(t,\sigma,v,w,u)|\leq
    k_1^*(t)+k_2^*(t)m+k_3^*(t)(|\sigma|_H+|v|_H+|w|_H)
    \end{equation*}  a.e. on $[0,T]$,  $\sigma,v,w,u\in H$.
\end{itemize}

Here, $k_i^*(t)=| k_i(t)|_{L^1(\Omega)}$, $i=1,2$, $k_j^*(t)=| k_j(t)|_{H}$, $j=0,3$.

\hspace{0.6cm} Let $\partial I_{\mathcal{K}(v,w)}(\sigma)$ be the subdifferential of
the indicator function of ${\mathcal{K}}(v,w)$ at a point $\sigma\in H$. Now we are in a position to define solutions for our control problems.

\begin{definition} A quadruple  $\{\sigma,v,w,u\}$ is called a solution of
control system $(\ref{1.1})$--$(\ref{1.6})$  if
\begin{itemize}
\item[$(i)$] $\sigma,v,w\in W^{1,2}(0,T;H)\cap
L^\infty(0,T;V)\cap L^2(0,T;H^2 (\Omega))$;

\item[$(ii)$] $u\in L^2(0,T;H)$;

\item[$(iii)$] $ \sigma'-(\lambda(v))'-\kappa
    \Delta_N\sigma+ \partial I_{\mathcal{K}(v,w)}(\sigma)
    \ni F(\sigma,v,w)\,u$ \; in $H$ a.e. on $[0,T]$;

\item[$(iv)$] $ v'- \Delta_N v=h(\sigma,v,w)$
\; in $H$ a.e. on $[0,T]$;

\item[$(v)$] $ w'-\Delta_N w=g(\sigma,v,w)$
\; in $H$ a.e. on $[0,T]$;

\item[$(vi)$] $ \sigma(0) = \sigma_0, \; v(0) = v_0, \; w(0) = w_0$ \; in  $H$;

\item[$(vii)$] $u(t)\in \mathcal{U}(t,v(t),w(t),\sigma(t))$ \; in $H$
a.e. on $[0,T]$,

\end{itemize}

where the prime denotes derivative with respect to $t$.

\hspace{0.6cm} A solution of control system $(\ref{1.1})$--$(\ref{1.5})$, $(\ref{1.7})$ is defined
similarly replacing the last inclusion with
\begin{equation*}
u(t)\in \overline{\co} \;\mathcal{U}(t,\sigma(t),v(t),w(t)) \quad \text{in} \;
H \;\; \text{for a.e.} \; t\in [0,T] .\label{}
\end{equation*}\end{definition}

\hspace{0.6cm}  The
sets of all solutions to control systems  $(\ref{1.1})$--$(\ref{1.6})$ and $(\ref{1.1})$--$(\ref{1.5})$, $(\ref{1.7})$ in the sense of Definition 2.1 we denote by $\mathcal{R}_\mathcal{U}$ and  $\mathcal{R}_{\overline{\co}\, \mathcal{U}}$, respectively.

\hspace{0.6cm} Defining the function
$\mathfrak{q}_{_\mathcal{U}}:[0,T]\times H^3\times H\to
{\mathbb{R}}\cup\{+\infty\}$  by the rule
\begin{equation*}\mathfrak{q}_{_\mathcal{U}}(t,\sigma,v,w,u)=
\left\{
\begin{array}{ll}
\mathfrak{q}(t,\sigma,v,w,u) & \mbox{if} \;u \in \mathcal{U}(t,\sigma,v,w) ,\\
+\infty  & \mbox{otherwise},
\end{array}\right.
\label{}\end{equation*} and denoting by
$\mathfrak{q}_{_\mathcal{U}}^{**}(t,\sigma,v,w,u)$ the bipolar of the function $u\mapsto
\mathfrak{q}_{_\mathcal{U}}(t,\sigma,v,w,u)$, one can easily see that
\begin{equation*}
\label{} \mathfrak{q}_{_\mathcal{U}}^{**}(t,\sigma,v,w,u)=\int_\Omega
q_{_{U}}^{**}(t,x,\sigma(x),v(x),w(x),u(x))\, dx.
\end{equation*}

So, our optimal control Problems $(P)$ and $(RP)$ can now be reformulated in the form:
\begin{itemize}
\item[(P)]  minimize \;\; $J(\sigma, v, w, u)={\displaystyle \int_0^T\mathfrak{q}(t,\sigma(t),v(t),w(t),u(t))\,dt}$ \;\; over  \,$\mathcal{R}$,

\item[(RP)]  minimize \;\; $J_{{U}}^{**}(\sigma, v, w, u)={\displaystyle \int_0^T\mathfrak{q}_{_\mathcal{U}}^{**}(t,\sigma(t),v(t),w(t),u(t))\,dt}$ \;\; over \, $\mathcal{R}_{\overline{\co}\, \mathcal{U}}$.

\end{itemize}

\hspace{0.6cm} Given Hypotheses $(H)$, $(U)$, and $(q)$, the main purpose of
this work is to prove the following results.

\begin{theorem} \label{Theorem2.1} For any  $(\sigma_*,v_*,w_*,u_*)\in \mathcal{R}_{\overline{\co}\,
\mathcal{U}}$ there exists a sequence $(\sigma_n,v_n,w_n,u_n)$ $\in\mathcal{R}_{\mathcal{U}}$, $n\geq 1$, such that
\begin{align}
(\sigma_n, v_n,w_n) &\to (\sigma_*,v_*,w_*) \quad \mbox{strongly in } C([0,T],H^3) ,\label{2.5}\\
  u_n &\to u_* \quad \mbox{weakly in } L^2([0,T],H) . \label{2.6}
\end{align}

Moreover,
\begin{align}
\nonumber \lim\limits_{n\rightarrow\infty}\sup\limits_{0\leq s\leq
t\leq T}
\Bigg|\int_{s}^{t}&\big(\mathfrak{q}_\mathcal{U}^{**}(\tau,\sigma_*(\tau),v_*(\tau),w_*(\tau),u_*(\tau))\\
&-\mathfrak{q}(\tau,\sigma_n(\tau),v_n(\tau),w_n(\tau),u_n(\tau))\big)\,
d\tau\Bigg| =0 .\label{2.7}
\end{align}
\end{theorem}

\begin{theorem} Problem \textup{(RP)} has an optimal solution and
\begin{equation} \min\limits_{(\sigma,v,w,u)\in \mathcal{R}_{\overline{\co} \,\mathcal{U}}} J^{**}_U(\sigma,v,w,u)
=\inf\limits_{(\sigma,v,w,u)\in \mathcal{R}_{\mathcal{U}}} J(\sigma,v,w,u) .\label{2.8}
\end{equation}
Moreover, for any solution $(\sigma_*,v_*,w_*,u_*)$ of  \textup{(RP)} there
exists a minimizing sequence $(\sigma_n,v_n,w_n,u_n)\in \mathcal{R}_{\mathcal{U}}$,
$n\geq 1$, for problem \textup{(P)}
such that \textup{(\ref{2.5})--(\ref{2.7})} hold.
\end{theorem}

\section{Properties of the control-to-state solution operator} \label{section3}

In this section, we define the control-to-state solution operator for our control systems and explore some of its properties which are crucial for establishing our main results in the next section. To this aim, first we let
$$S_m:=\{u\in L^2(0,T;H); \;  |u(t,x)|  \leq m \quad \mbox{a.e. on} \; Q(T)\}.
$$

Due to convexity, the bound for the control functions of Problem $(P)$ given in  Hypothesis $(U3)$   obviously extends to the control functions  of the convexified problem $(RP)$. In particular,  the controls of both problems
belong to the set $S_m$. Accordingly, let  $\mathcal{T}:S_m\to C([0,T];H\times H\times H)$ be the operator
which with each $u\in S_m$ associates the unique solution of
system $(\ref{1.1})$--$(\ref{1.5})$:
\begin{equation}
(\sigma_u,v_u,w_u)=\mathcal{T}(u). \label{3.1}
\end{equation}

The existence and uniqueness of such a solution as well as uniform a priori
estimates for all possible solutions  independent of the
control $u$ are provided by the following theorem.

\begin{theorem} \label{th3-1}
 For any fixed $u\in S_m$ system
$(\ref{1.1})$--$(\ref{1.5})$ has a unique solution. Moreover, for
any solution $(v,w,\sigma)$ of $(\ref{1.1})$--$(\ref{1.5})$ with $u\in S_m$ the following estimates
\begin{equation} \quad 0 \leq \sigma,v, w
\leq M_0 \quad \mbox{ a.e. on } Q(T), \label{3.2}
\end{equation}
\begin{align}\nonumber & |\sigma' |_{L^2(0,T;
H)}+ | v' |_{L^2(0,T; H)} + |w' |_{L^2(0,T;
H)}\\ \nonumber  + & |\Delta \sigma|_{L^2(0,T; H)} + |\Delta v|_{L^2(0,T; H)}+ |\Delta w|_{L^2(0,T; H)}   \\   + &| \nabla \sigma|_{L^{\infty}(0,T; H)} + | \nabla v|_{L^{\infty}(0,T; H)}+ | \nabla w|_{L^{\infty}(0,T; H)}
\leq M_0 \label{3.3}
\end{align}
hold for a  constant $M_0>0$ independent of $u$.
\end{theorem}

\bpf{Proof} The existence of a unique solution to
$(\ref{1.1})$--$(\ref{1.5})$ with a fixed $u\in S_m$ and the estimate (\ref{3.2})  can be proved
following the pattern  of \cite[Theorems 3.1, 3.2, and 3.10]{AM}.

\hspace{0.6cm} By virtue of the bound (\ref{3.2}),
 we may now assume (cutting off outside
the set where $\sigma$, $v$, and $w$ are bounded, if necessary) that the
functions $ F, h, g$ are all bounded  (with a common bound $M>0$) and globally Lipschitz continuous $($with a common Lipschitz constant $L>1)$.

\hspace{0.6cm} Below, we recap a part of the reasoning of \cite{AM} which will allow us to establish the uniform energy estimates (\ref{3.3}). To this end, given $\mu>0$ we introduce the following approximate system:
\begin{equation}
\*\hspace{-0.1cm}\sigma_\mu'-(\lambda(v_\mu))'-\kappa
    \Delta_N\sigma_\mu+ \partial I^\mu_{\mathcal{K}(v_\mu,w_\mu)}(\sigma_\mu)
    = F(\sigma_\mu,v_\mu,w_\mu)\,u \;\, \text{in} \; H \; \text{a.e. on} \; [0,T]
,\label{3.4}
\end{equation}
\begin{equation} v_\mu'- \Delta_N v_\mu=h(\sigma_\mu,v_\mu,w_\mu) \qquad \text{in} \; H \; \text{a.e. on} \; [0,T] , \label{3.5}\end{equation}
\begin{equation} w_\mu'-\Delta_N w_\mu=g(\sigma,v_\mu,w_\mu)  \qquad \text{in} \; H \; \text{a.e. on} \; [0,T] . \label{3.6}\end{equation}

Here, for $\mu>0$ and $\sigma, v, w\in \mathbb{R}$ the function $\partial I^\mu_{\mathcal{K}(v,w)}(\sigma)$ is the Yosida regularization of the subdifferential $\partial I_{\mathcal{K}(v,w)}(\sigma)$.

\hspace{0.6cm} By \cite{AM}, for any $\mu>0$ there exists a unique triplet $(\sigma_\mu, v_\mu, w_\mu)$ solving the approximate system above and satisfying the initial conditions $\sigma_\mu(0)=\sigma_0$, $v_\mu(0)=v_0$, $w_\mu(0)=w_0$. Arguing similarly to \cite[3.3.2. Proof of Theorem 3.6]{{AM}} we obtain the following counterparts of the inequalities $(21), (22), (23), (32)$ of this reference for such triplets $(\sigma_\mu, v_\mu, w_\mu)$, $\mu>0$:
\begin{equation} |v_\mu'|^2_H +|w_\mu'|^2_H +\frac{d}{dt}|\nabla v_\mu|^2_H+\frac{d}{dt}|\nabla w_\mu|^2_H\leq C_1 ,  \label{21}\end{equation}

\begin{equation} \frac{d}{dt}|\nabla v_\mu|^2_H+\frac{d}{dt}|\nabla w_\mu|^2_H+|\Delta v_\mu|^2_H +|\Delta w_\mu|^2_H \leq C_1 ,  \label{22}\end{equation}

\begin{equation} |\sigma_\mu'|^2_H + \kappa\frac{d}{dt}|\nabla \sigma_\mu|^2_H+2\frac{d}{dt}I^\mu_{\mathcal{K}(v_\mu,w_\mu)}(\sigma)\leq C_2 (|v_\mu'|^2_H +|w_\mu'|^2_H+\kappa^2|\Delta \sigma_\mu|^2_H +1) ,  \label{23}\end{equation}

\begin{align}\nonumber \frac{d}{dt}\bigg\{I^\mu_{v_\mu,w_\mu}(\sigma_\mu)&+\frac{1}{2}|\nabla \sigma_\mu|^2_H-(\lambda'(v_\mu)\nabla v_\mu,\nabla \sigma_\mu)_H\bigg\} \\ \nonumber & +\kappa|\Delta \sigma_\mu|^2_H +\frac{1-\kappa}{4}|\partial I^\mu_{\mathcal{K}(v_\mu,w_\mu)}(\sigma_\mu)|^2_H \\ \nonumber &\leq C_3(1+|\nabla \sigma_\mu|^2_H+|\nabla v_\mu|^2_H+|\nabla w_\mu|^2_H+|\sigma_\mu'|^2_H+|v_\mu'|^2_H +|w_\mu'|^2_H\\
&\hspace{1.21cm} +(1+\kappa)(1+|\Delta v_\mu|^2_H +|\Delta w_\mu|^2_H )).  \label{32}\end{align}

a.e. on $(0,T)$, where $C_1, C_2, C_3$  are some positive constants depending on  $|\lambda'|_{L^\infty(\mathbb{R})}$, $|\lambda''|_{L^\infty(\mathbb{R})}$, $|f_*|_{W^{2,\infty}(\mathbb{R}^2)}$, $|f_*|_{W^{2,\infty}(\mathbb{R}^2)}$, $L$,  $M_0$, $M$, $m$, and $\mu(\Omega)$, but independent of $\mu$. Next, calculating $(\ref{21})+\varepsilon_1\times(\ref{22})+\varepsilon_2\times(\ref{23})+\varepsilon_3\times(\ref{32})$ with positive numbers $\varepsilon_1, \varepsilon_2, \varepsilon_3$
to be specified later we obtain
\begin{align}\nonumber & (\varepsilon_2 C_2-\varepsilon_3 C_3) |\sigma_\mu'|^2_H+(1-\varepsilon_2 C_2-\varepsilon_3 C_3)|v_\mu'|^2_H +(1-\varepsilon_2 C_2-\varepsilon_3 C_3) |w_\mu'|^2_H\\
\nonumber & +\kappa(\varepsilon_3-\varepsilon_2 C_2\kappa)|\Delta \sigma_\mu|^2_H+(\varepsilon_1-\varepsilon_3 C_3(1+\kappa))|\Delta v_\mu|^2_H +(\varepsilon_1-\varepsilon_3 C_3(1+\kappa))|\Delta w_\mu|^2_H  \\
\nonumber & +\varepsilon_3 \frac{1-\kappa}{4}|\partial I^\mu_{\mathcal{K}(v_\mu,w_\mu)}(\sigma_\mu)|^2_H +\frac{d}{dt}\bigg\{(2\varepsilon_2+\varepsilon_3) I^\mu_{\mathcal{K}(v_\mu,w_\mu)}(\sigma_\mu)-\varepsilon_3(\lambda'(v_\mu)\nabla v_\mu,\nabla \sigma_\mu)_H\\
\nonumber &+(\varepsilon_2 \kappa+\frac{\varepsilon_3}{2})|\nabla \sigma_\mu|^2_H+(1+\varepsilon_1)|\nabla v_\mu|^2_H+(1+\varepsilon_1)|\nabla w_\mu|^2_H \bigg\}\\
\nonumber & \leq C_4 +\varepsilon_3 C_3(2+\kappa+|\nabla \sigma_\mu|^2_H+|\nabla v_\mu|^2_H+|\nabla w_\mu|^2_H),
 \label{}\end{align}

where $C_4=C_1+\varepsilon_1 C_1+\varepsilon_2 C_2$. Integrating this inequality from $0$ to $t\in (0,T]$ and then estimating  the term $-\varepsilon_3(\lambda'(v_\mu)\nabla v_\mu,\nabla \sigma_\mu)_H$ in the resulting inequality as follows
$$-\varepsilon_3(\lambda'(v_\mu)\nabla v_\mu,\nabla \sigma_\mu)_H\geq -\frac{1}{4}\varepsilon_3^2|\lambda'|_{L^\infty(\mathbb{R})}^2|\nabla \sigma_\mu|_H^2-|\nabla v_\mu|_H^2$$

we infer that
\begin{align}\nonumber (\varepsilon_2 C_2-\varepsilon_3 C_3) |\sigma_\mu'|^2_{L^2(0,t;H)}&+(1-\varepsilon_2 C_2-\varepsilon_3 C_3)|v_\mu'|^2_{L^2(0,t;H)} \\
\nonumber &+(1-\varepsilon_2 C_2-\varepsilon_3 C_3) |w_\mu'|^2_{L^2(0,t;H)} \\
\nonumber   +\kappa(\varepsilon_3-\varepsilon_2 C_2\kappa)|\Delta \sigma_\mu|^2_{L^2(0,t;H)}&+(\varepsilon_1-\varepsilon_3 C_3(1+\kappa))|\Delta v_\mu|^2_{L^2(0,t;H)} \\
\nonumber & +(\varepsilon_1-\varepsilon_3 C_3(1+\kappa))|\Delta w_\mu|^2_{L^2(0,t;H)} \\
\nonumber & \hspace{-4.62cm} +\varepsilon_3 \frac{1-\kappa}{4}|\partial I^\mu_{\mathcal{K}(v_\mu,w_\mu)}(\sigma_\mu)|^2_{L^2(0,t;H)} +(2\varepsilon_2+\varepsilon_3) I^\mu_{\mathcal{K}(v_\mu,w_\mu)}(\sigma_\mu)\\
\nonumber & \hspace{-4.62cm} +(\varepsilon_2 \kappa+\frac{\varepsilon_3}{2}-\frac{1}{4}\varepsilon_3^2|\lambda'|_{L^\infty(\mathbb{R})}^2)|\nabla \sigma_\mu|^2_H+\varepsilon_1|\nabla v_\mu|^2_H+(1+\varepsilon_1)|\nabla w_\mu|^2_H \\
 & \hspace{-4.62cm} \leq C_5 +\varepsilon_3 C_3 \int_0^t \Big\{|\nabla \sigma_\mu(\tau)|^2_H+|\nabla v_\mu(\tau)|^2_H+|\nabla w_\mu(\tau)|^2_H\Big\} d\tau,
 \label{Ineq}\end{align}

where $C_5$ is a positive constant which depends on $|\sigma_0|_V$, $|v_0|_V$, $|w_0|_V$, but is independent of $\mu$. Now, all the coefficients in $(\ref{Ineq})$ will be positive provided we take
$\varepsilon_1=\min\{\min\{1, \frac{1}{2C_2}\},  \frac{4C_3}{|\lambda'|_{L^\infty(\mathbb{R})}^2}\}$, $\varepsilon_2= \frac{1}{4C_2}$, $\varepsilon_3=\min\{\frac{1}{4C_3}\min\{1, \frac{1}{2C_2}\},  \frac{1}{|\lambda'|_{L^\infty(\mathbb{R})}^2}\}$ and
choose $\kappa\in (0,\kappa_0)$ with $\kappa_0=\min\{\frac{1}{2},\frac{1}{2C_3}\min\{1, \frac{1}{2C_2}\},  \frac{2}{|\lambda'|_{L^\infty(\mathbb{R})}^2}\}$. Therefore, invoking Gronwall's inequality from $(\ref{Ineq})$ we obtain the following uniform energy estimates for the triplets $(\sigma_\mu, v_\mu, w_\mu)$, $\mu>0$, solving the approximate problem (\ref{3.4})--(\ref{3.6}):
\begin{align}\nonumber & |\sigma_\mu' |_{L^2(0,T;
H)}+ | v'_\mu |_{L^2(0,T; H)} + |w'_\mu |_{L^2(0,T;
H)}\\ \nonumber  + & |\Delta \sigma_\mu|_{L^2(0,T; H)} + |\Delta v_\mu|_{L^2(0,T; H)}+ |\Delta w_\mu|_{L^2(0,T; H)}   \\   \nonumber + &| \nabla \sigma_\mu|_{L^{\infty}(0,T; H)} + | \nabla v_\mu|_{L^{\infty}(0,T; H)}+ | \nabla w_\mu|_{L^{\infty}(0,T; H)} \\ + & |\partial I^\mu_{\mathcal{K}(v_\mu,w_\mu)}(\sigma_\mu)|_{L^2(0,T;H)}+| I^\mu_{\mathcal{K}(v_\mu,w_\mu)}(\sigma_\mu)|_{L^\infty(0,T)}
\leq M_0, \label{3.12}
\end{align}

By the weak and weak-star compactness
results, this bound allows us to conclude  that there exists a  null sequence $\mu_n$, $n\geq 1$,  and functions $\sigma,v,w\in W^{1,2}(0,T;H)\cap
L^\infty(0,T;V)\cap L^2(0,T;H^2(\Omega))$  such
that
\begin{align}
\sigma_n:=\sigma_{\mu_{n}} \to \sigma,  \quad v_n:=v_{\mu_{n}} \to v,  \quad w_n:=w_{\mu_{n}} \to w \label{3.13}
\end{align}

weakly in $W^{1,2}(0,T; H) \cap
L^2(0,T;H^2(\Omega))$ and weakly-star in $L^{\infty}(0,T; V)$, and, thus, strongly in $C([0,T];H)$.

\hspace{0.6cm} With the convergences $(\ref{3.13})$ at hand, we can now pass to the limit in Eqs. $(\ref{3.5})$, $(\ref{3.6})$ to infer that the triplet $(\sigma,v,w)$ satisfies Eqs. $(iv)$, $(v)$ of Definition 2.1. Furthermore, denoting $f_n:=-\sigma_n'+(\lambda(v_n))'+\kappa\Delta_N\sigma_n+ F(\sigma_n,v_n,w_n)\,u$ we see that
\begin{align}
f_n\to f:=-\sigma'+(\lambda(v))'+\kappa\Delta_N\sigma+ F(\sigma,v,w)\,u  \label{3.14}
\end{align}

weakly in $L^{2}(0,T;H)$. Since \cite[Proposition 0.3.5]{Kenmochi}
$$\partial I^{\mu_n}_{\mathcal{K}(v_n,w_n)}(\sigma_n)\in \partial I_{\mathcal{K}(v_n,w_n)}(J_{\mathcal{K}(v_n,w_n)}(\sigma_n)),$$

where $J_{\mathcal{K}(v_n,w_n)}(\sigma_n)$ is the projection of $\sigma_n$ onto the set ${\mathcal{K}(v_n,w_n)}$, we also have
\begin{align}
f_n\in \partial I_{\mathcal{K}(v_n,w_n)}(J_{\mathcal{K}(v_n,w_n)}(\sigma_n)).  \label{3.15}
\end{align}

On the other hand, the boundedness of $I^{\mu_n}_{\mathcal{K}(v_n,w_n)}(\sigma_n)$ implies that
\begin{align}
J_{\mathcal{K}(v_n,w_n)}(\sigma_n)\to \sigma  \label{3.16}
\end{align}

strongly in $C([0,T];H)$. In view of (\ref{3.13})--(\ref{3.16}), the application of Lemma 2.1 yields
\begin{align*}
f\in \partial I_{\mathcal{K}(v,w)}(\sigma),  \label{}
\end{align*}

so that the triplet $(\sigma,v,w)$ satisfies Eqs. $(iii)$ of Definition 2.1 as well. We note that inequality (\ref{2.1}) in this case follows from the Lipschitzness of $f_*$ and $f^*$. Therefore, we conclude that $(\sigma,v,w)$ is a solution of system (\ref{1.1})--(\ref{1.5}) with a fixed $u\in S_m$. The estimate (\ref{3.3}) now follows from (\ref{3.12}) and (\ref{3.13}).
\epf

\begin{theorem} \label{Theorem 3.2}
The operator $\mathcal{T}:S_m\to C([0,T];H\times H\times H)$
is weak-strong continuous.
\end{theorem}

\bpf{Proof} We note that the set $S_m$ when endowed with the weak topology of
the space $L^2(0,T;H)$ is metrizable. Consequently, to prove the theorem it is enough to show the
sequential continuity of the operator $\mathcal{T}$. Let, then,    $u_n$,
$n\geq 1$, be  an arbitrary sequence from $S_m$  weakly converging  to
some $u\in S_m$. Denote by $(\sigma_n,v_n,w_n):=(\sigma(u_n),v(u_n),w(u_n))$, $n\geq 1$,
the sequences of solutions to system  $(\ref{1.1})$--$(\ref{1.5})$
corresponding to the controls $u_n$, $n\geq 1$. Similarly as in the proof of Theorem 3.1, by the weak and weak-star compactness
results, the uniform estimates (\ref{3.2}) and (\ref{3.3}) imply that there exists a subsequence $(\sigma_k,v_k,w_k):=(\sigma_{n_k},v_{n_k},w_{n_k})$, $k\geq 1$, of
the sequence $(\sigma_n,v_n,w_n)$, $n\geq 1$, and some functions $\sigma,v,w\in W^{1,2}(0,T;H)\cap
L^\infty(0,T;V)\cap L^2(0,T;H^2(\Omega))$  such
that
\begin{align}
\sigma_k \to \sigma,  \qquad v_k \to v,  \qquad w_k \to w \label{3.17}
\end{align}

weakly in $W^{1,2}(0,T; H) \cap
L^2(0,T;H^2(\Omega))$ and weakly-star in $L^{\infty}(0,T; V)$, and, thus, strongly in $C([0,T];H)$.

\hspace{0.6cm} From the Lipschitz continuity of the functions $F,g,h$ and $(\ref{3.17})$ we deduce that
\begin{align}\nonumber -\sigma_k'+ (\lambda(v_k))'+\kappa\Delta_N\sigma_k &+ F(\sigma_k,v_k,w_k)u_k \\ &\to -\sigma'+ (\lambda(v))'+\kappa\Delta_N\sigma+ F(\sigma,v,w)u
\label{3.18}
\end{align}

weakly in $L^{2}([0,T]; H)$ and
\begin{align}
 h(\sigma_k,v_k,w_k)\to h(\sigma,v,w),\quad g(\sigma_k,v_k,w_k)\to g(\sigma,v,w)
 \label{3.19}
\end{align}

strongly in $C([0,T];H)$. The claim of the theorem now follows from the convergences (\ref{3.17})--(\ref{3.19}) and Lemma 2.1.
\epf

\begin{theorem} \textup{\cite[Theorem 3.2]{SCL2019}} \label{Theorem3.1}
Let $u_i \in S_m$  and $\{\sigma_i,v_i, w_i\} = \mathcal{T}(u_i)$, $i = 1, 2$. Then,
\begin{align}
 |\sigma_1(t)-\sigma_2(t)|_H^2+ |v_1(t) - v_2(t)|_H^2  &+ |w_1(t)-w_2(t)|_H^2  \nonumber \\
&\leq  \  C_m   \int_0^t |u_1(\tau)-u_2(\tau)|_H^2\,d\tau,  \label{3.20}
\end{align}

$t\in [0,T]$, where $C_m>0$ is a positive constant depending only on $m$.
\end{theorem}

\section{Proofs of the main results}

In this section, we prove Theorems 2.1 and 2.2. To this end, first, on the basis of the control constraint multivalued mapping
$\mathcal{U}$ and the cost integrand $\mathfrak{q}$ we construct an auxiliary multivalued mapping $\mathcal{U}_{\mathfrak{q}}:[0,T]\times H^3\rightarrow
{H}\times \mathbb{\mathbb{R}}$ as follows
\begin{equation}\mathcal{U}_{\mathfrak{q}}(t,\sigma,v,w)=\{(u,\mathfrak{q}(t,\sigma,v,w,u))\in
{H}\times \mathbb{R};\; u\in \mathcal{U}(t,\sigma,v,w)\},
\label{4.1}\end{equation}

where the Banach space
$H\times\mathbb{R}$ is equipped with the norm
\begin{equation}
\label{4.2}|(u,\lambda)|_{{H}\times \mathbb{R}}:=\max\{|u|_H,|\lambda|\},\qquad
u\in H, \lambda\in \mathbb{R}.\end{equation}

Given the properties  $(\mathcal{U}\textbf{1})$--$(\mathcal{U}\textbf{3})$ and $(\mathfrak{q}\textbf{1})$--$(\mathfrak{q}\textbf{3})$ we can easily show that the  mapping ${\mathcal{U}}_\mathfrak{q}$ is
measurable in $t$, continuous in the Hausdorff metric on the space ${\rm
cb}(H\times \mathbb{R})$ in $(\sigma,v,w)$ a.e. on $[0,T]$ and
\begin{align}\nonumber|{\mathcal{U}}_\mathfrak{q}(t,\sigma,v,w)&|_{{H}\times \mathbb{R}}\\ &\leq m\mu(\Omega)^{1/2}+k_1^*(t)+k_2^*(t)m+k_3^*(t)(|\sigma|_H+|v|_H+|w|_H).  \label{4.3}
\end{align}

Moreover, we have (cf. \cite[Lemma 2.2]{Tolstonogov2004})

\begin{equation}\*\hspace{-0.4cm}\mathfrak{q}_\mathcal{U}^{**}(t,\sigma,v,w,u)=\*\hspace{-0.1cm} \left\{
\begin{array}{ll}
\*\hspace{-0.2cm}\min\{\lambda\in \mathbb{R}; (u,\lambda)\in \overline{\co}\,{\mathcal{U}}_\mathfrak{q}(t,\sigma,v,w) \} & \*\hspace{-0.3cm}\mbox{if} \;u \in \overline{\co}\,{\mathcal{U}}(t,\sigma,v,w) ,\\
\*\hspace{-0.2cm}+\infty  & \mbox{otherwise},
\end{array}\right.
\label{4.4}\end{equation}

and for any $\varepsilon>0$, there is  a closed set $T_\varepsilon\subset [0,T]$, $\mu([0,T]\setminus T_\varepsilon)\leq\varepsilon$, such
that  $\mathfrak{q}_\mathcal{U}^{**}$ restricted to $T_\varepsilon\times H^3\times H$ is lower semicontinuous. This implies, in particular, that given \emph{an arbitrary} $(\sigma_*,v_*,w_*,u_*)\in
\mathcal{R}_{\overline{\co}\, \mathcal{U}}$, the function $t\to
\mathfrak{q}_\mathcal{U}^{**}(t,\sigma_*(t),v_*(t),w_*(t),u_*(t))$ is
measurable and
$$
\left(u_*(t),\mathfrak{q}_\mathcal{U}^{**}(t,\sigma_*(t),v_*(t),w_*(t),u_*(t))\right)\in
\overline{\co}\, {\mathcal{U}}_\mathfrak{q}(t,\sigma_*(t),v_*(t),w_*(t))  $$
a.e. on $[0,T]$. Since the bound in (\ref{4.3}) evidently extends to $\overline{\co}\, {\mathcal{U}}_\mathfrak{q}$, invoking \cite[Corollary 1.1]{Chuong}, from the last inclusion we obtain the existence of a measurable function
${\gamma}_n(t)\in {\mathcal{U}}(t,\sigma_*(t),v_*(t),w_*(t))
$
which, in view of (\ref{4.2}), satisfies
\begin{equation}
\sup\limits_{0\leq s\leq t\leq T} \left|\int_s^{t}
\left({u}_*(\tau)-{\gamma}_n(\tau)\right) \, d\tau\right|_H \leq
\frac{1}{n} ,\label{4.5}
\end{equation}
and
\begin{align}\nonumber \sup\limits_{0\leq s\leq t\leq T}
\Bigg|\int_s^{t}
&\big(\mathfrak{q}_\mathcal{U}^{**}(\tau,\sigma_*(t),v_*(\tau),w_*(\tau),u_*(\tau))\\
&-\mathfrak{q}(\tau,\sigma_*(t),v_*(\tau),w_*(\tau),{\gamma}_n(\tau))\big) \,
d\tau\Bigg| \leq \frac{1}{n} .\label{4.6}
\end{align}

\hspace{0.6cm}  From Theorem 3.2 it follows that the set
\begin{equation*} \mathcal{R}:=\{(\sigma,v,w)\in C([0,T];H^3); \; (\sigma,v,w)=\mathcal{T}(u), \; u\in S_m\}\label{}\end{equation*}

is compact in $ C([0,T];H^3)$. Fix $n\geq 1$. The property  $(\mathcal{U}2)$ implies that for any $(\sigma,v,w)\in H^3$ and a.e. $t\in [0,T]$ there exists ${\gamma}\in {\mathcal{U}}(t,\sigma,v,w)$ such that
\begin{align*}|\gamma_n(t)-\gamma|_H^2&< \frac{2}{n^2}+2 d^{\,2}_H(\gamma_n(t),{\mathcal{U}}(t,\sigma,v,w))\\
&< \frac{2}{n^2}+12k^2(t)\left(|\sigma_*(t)-\sigma|_H^2+|v_*(t)-v|_H^2+|w_*(t)-w|_H^2\right). \label{}\end{align*}

Define the multivalued mapping
\begin{align}\nonumber \mathcal{U}_n(t,\sigma,v,w)&:=\{ \gamma\in \mathcal{U}(t,\sigma,v,w); \; |\gamma_n(t)-\gamma|_H^2 \\
&\leq \frac{2}{n^2}+12k^2(t)\left(|\sigma_*(t)-\sigma|_H^2+|v_*(t)-v|_H^2+|w_*(t)-w|_H^2\right)\}, \label{4.7}\end{align}

and the associated Nemytskii multivalued operator $\Gamma_n:\mathcal{R}\to L^2(0,T;H)$:
\begin{equation} \*\hspace{-0.4cm}\Gamma_n(\sigma,v,w):=\{u\in L^2(0,T;H); \; u(t)\in {\mathcal{U}_n(t,\sigma(t),v(t),w(t))} \; \mbox{for a.e.} \; t\in[0,T]\},   \label{4.8}\end{equation}

which is lower semicontinuous with nonempty closed decomposable
values. According to \cite[Theorem 3.1]{Fryszkowski} there exists a continuous mapping
${\alpha}_n:\mathcal{R}\to L^1(0,T;H)$ such that
\begin{equation} {\alpha}_n(\sigma,v,w)\in \Gamma_n(\sigma,v,w), \; (\sigma,v,w)\in \mathcal{R} .\label{4.9}\end{equation}

In view of the property $(\mathcal{U}3)$ and the definition of $\mathcal{U}_n$, from this inclusion we deduce that
 ${\alpha}_n$ is continuous from $\mathcal{R}$ to $L^2(0,T;H)$ as well
and ${\alpha}_n(\sigma,v,w) \in S_m$, $(\sigma,v,w)\in \mathcal{R}$.

\hspace{0.6cm} Consider now the superposition of $\mathcal{T}$ and
${\alpha}_n$. By virtue of Theorem 3.2 it follows that this superposition
${\alpha}_n\circ \mathcal{T}:S_m\to S_m$ is weak-weak continuous.
Since $S_m$ is obviously convex and compact in the weak topology of
the space $L^2(0,T;H)$, from the Schauder fixed point theorem it follows that there exists a fixed point $u_n\in S_m$ of the
operator ${\alpha}_n\circ \mathcal{T}$, i.e.
\begin{equation} u_n={\alpha}_n(\mathcal{T}(u_n)) .\label{4.10}\end{equation}

Setting $(\sigma_n,v_n,w_n):=\mathcal{T}(u_n)$, from (\ref{4.7})--(\ref{4.10})
we see that $(\sigma_n,v_n,w_n,u_n)\in \mathcal{R}_\mathcal{U}$, $n\geq 1$,  and
\begin{align}\nonumber
|&{\gamma}_n(t)-u_n(t)|_H^2 \\ &\leq \frac{2}{n^2}
+12k^2(t)\left(|\sigma_*(t)-\sigma_n(t)|_H^2+|v_*(t)-v_n(t)|_H^2+|w_*(t)-w_n(t)|_H^2\right).
\label{4.11}
\end{align}

Since on the set $S_m$ the topology induced by the ``weak'' norm in (\ref{4.5}) coincides with the weak topology of the space $L^2(0,T;H)$,
from (\ref{4.5}) we obtain
\begin{align}
\gamma_n \to u_* \quad \mbox{weakly in } L^2(0,T;H)  .\label{4.12}
\end{align}

Theorem 3.2 implies then that
\begin{align}(\sigma(\gamma_n),v(\gamma_n),w(\gamma_n)) \to (\sigma_*,v_*,w_*) \quad \mbox{strongly in }
C([0,T];H^3) . \label{4.13}
\end{align}

From (\ref{4.11}) and (\ref{3.20}) it follows that
\begin{align*}
|\sigma_n(t)&-\sigma(\gamma_n)(t)|_H^2 +|v_n(t)-v(\gamma_n)(t)|_H^2 +|w_n(t)-w(\gamma_n)(t)|_H^2\\ &\leq
C_m\int_0^t
|u_n(\tau)-\gamma_n(\tau)|_H^2 d\tau \\
 &\leq C_m\int_0^t
\bigg(\frac{2}{n^2}+24k^2(\tau)\big(|\sigma_n(\tau)-\sigma(\gamma_n)(\tau)|_H^2
+|v_n(\tau)-v(\gamma_n)(\tau)|_H^2 \\
&\*\hspace{1.9cm}+ |w_n(\tau)-w(\gamma_n)(\tau)|_H^2\big)
+24k^2(\tau)\big(|\sigma(\gamma_n)(\tau)-\sigma_*(\tau)|_H^2 \\
&\*\hspace{1.9cm}+|v(\gamma_n)(\tau)-v_*(\tau)|_H^2 +
|w(\gamma_n)(\tau)-w_*(\tau)|_H^2\big)\bigg) d\tau.\end{align*}

Invoking Gronwall's  lemma,  from the last inequality and
(\ref{4.13}) we obtain
\begin{equation*}
(\sigma_n,v_n,w_n) \to (\sigma_*,v_*,w_*) \quad \mbox{strongly in } C([0,T];H^3)
.
\end{equation*}

Finally, from (\ref{4.11}), (\ref{4.12}), and Lebesgue's dominated
convergence theorem we conclude that
\begin{equation*}
u_n \to u_* \quad \mbox{weakly in } L^2(0,T;H),
\end{equation*}

and the claim of Theorem 2.1 follows.

\hspace{0.6cm} Set now
$\mathcal{R}:=\{(\sigma(t),v(t),w(t));\; t\in [0,T], \, (\sigma,v,w)\in \mathcal{R}_\mathcal{U}\}.$
Theorem 3.2 together with Hypothesis $(U3)$ imply that the set
$\mathcal{R}$ is compact in $H^3$.
Define $\varphi:[0,T]\times
H^3\times H\to {\mathbb{R}}\cup\{+\infty\}$ by
\begin{equation*}\varphi(t,\sigma,v,w,u)= \left\{
\begin{array}{ll}
\mathfrak{g}_\mathcal{U}^{**}(t,\sigma,v,w,u), & \mbox{if }t\in [0,T], \;  (\sigma,v,w)\in \mathcal{R}, \; u \in H,\\
+\infty , & \mbox{otherwise}.
\end{array}\right.
\label{}\end{equation*}

By the property of $\mathfrak{g}_\mathcal{U}^{**}$ stated after (\ref{4.4}), there exists a
 sequence of  closed sets
$T_n\subset T_{n+1}\subset \dots\subset[0,T]$, $n\geq 1$, with $\mu([0,T]\setminus
\bigcup_{n=1}^\infty T_n)=0$ such that the function
$\mathfrak{g}_\mathcal{U}^{**}(t,v,w,u)$ restricted to $T_n\times
\mathcal{R}\times H$, $n\geq 1$, is lower semicontinuous.
Since the set
$\bigcup_{n=1}^\infty T_n$ is a Borel, the
function $\varphi$ is Borel on $\bigcup_{n=1}^\infty T_n\times
\mathcal{R}\times H$. Letting the values of $\varphi$ to be zero outside the set  $\bigcup_{n=1}^\infty T_n\times
\mathcal{R}\times H$, we can assume that $\varphi$ is Borel on  $[0,T]\times
\mathcal{R}\times H$ and, thus, it is measurable.  From the fact that $\mathfrak{g}_\mathcal{U}^{**}$ is the largest convex function such that $\mathfrak{g}_\mathcal{U}^{**}\leq \mathfrak{g}$ and from
$(\mathcal{U}\textbf{3})$, $(\mathfrak{q}\textbf{3})$, (\ref{4.4}) it follows that
$$-l(t)\leq q(t,\sigma,v,w,u) \quad\text{for a.e.} \;t\in [0,T], \; \sigma,v,w,u\in H,$$

for some $l\in L^1([0,T],\mathbb{R}^+)$.

\hspace{0.6cm} Defining now the
integral functional
$$(v,w,u)\to J_\varphi(\sigma,v,w,u):=\int_0^T \varphi(t,\sigma(t),v(t),w(t),u(t))\, dt$$
from  \cite[Theorem 2.1]{Balder} we deduce that $J_\varphi$ is sequentially lower semicontinuous on the space $C([0,T],H^3)\times L^{2}([0,T],H)$, when $L^{2}([0,T],H)$ is endowed with the weak topology.
According to Theorem
3.2 the set $\mathcal{R}_{\overline{\co} \,\mathcal{U}}$ is
compact in this space.
Since, obviously,
$J_\varphi(\sigma,v,w,u)=J^{**}_U(\sigma,v,w,u)$ on
$\mathcal{R}_{\overline{\co} \, \mathcal{U}}$, we conclude that
problem $(RP)$ has a solution $(\sigma_*,v_*,w_*,u_*)\in
\mathcal{R}_{\overline{\co} \,\mathcal{U}}$. Then,  Theorem
2.1 implies that  there exists a sequence $(\sigma_n,v_n,w_n,u_n)\in
\mathcal{R}_{\mathcal{U}}$, $n\geq 1$, such that
(\ref{2.5})--(\ref{2.7}) hold. In particular,
\begin{equation*}
J^{**}_U(\sigma_*,v_*,w_*,u_*)=\lim\limits_{n\to\infty}J(\sigma_n,v_n,w_n,u_n),
\end{equation*}
which combined with
\begin{equation*}
J^{**}_U(\sigma_*,v_*,w_*,u_*)\leq J^{**}_U(\sigma_n,v_n,w_n,u_n) \leq J(\sigma_n,v_n,w_n,u_n)
\end{equation*}
proves Theorem 2.2.

\bibliographystyle{amsplain}

\end{document}